\documentclass[a4paper,12pt]{article}
\usepackage{amsfonts}
\usepackage{amsmath, amsfonts, amssymb, amsthm, amscd}
\usepackage{index}
\usepackage[T1]{fontenc}
\usepackage{t1enc}

\usepackage{bbm}
\usepackage{graphicx}
\usepackage{fancybox}
\usepackage{fancyhdr}%
\usepackage{empheq}
\usepackage{bm}
\usepackage{xcolor}
\usepackage{titlesec}
       \titleformat{\chapter}[display]
             {\normalfont\Large\bfseries}{\thechapter}{12pt}{\Large}
       \titleformat{\section}
             {\normalfont\large\bfseries}{\thesection}{12pt}{\large}
       \titlespacing*{\chapter}{0pt}{0pt}{15pt} 
       \titlespacing*{\section}{0pt}{3.5ex plus 1ex minus .2ex}{2.3ex plus .2ex}

\usepackage[titletoc]{appendix}

\renewcommand{\d}{\textrm{d}}

\newcommand{\vth}{\vartheta}
\newcommand{\vep}{\varepsilon}

\newcommand{\E}{\textbf{\textrm{E}}}
\renewcommand{\P}{\textbf{{\textrm{P}}}}
\def\1{\mbox{1\hspace{-.25em}I}}
\newtheorem{theorem}{Theorem}[section]

\newtheorem{lemma}[theorem]{Lemma}
\newtheorem{prop}[theorem]{Proposition}

\numberwithin{equation}{section}

\begin{document}
\date{}
\title{On Goodness-of-fit Testing for Ergodic Diffusion Process with
  Shift Parameter\footnote{This work has been partially supported by MIUR grant 2009.}}
\author{Ilia Negri, Li
Zhou\footnote{Corresponding
    author. \textit{E-mail address}:
li.zhou.etu@univ-lemans.fr. \newline\newline Ilia Negri\newline
Department of Information Technology and Mathematical Methods,
University of Bergamo. Viale Marconi 5, 24044 Dalmine (BG)
Italy.\newline\newline Li Zhou\newline Laboratoire de Statistique et
Processus, Universit\'e du Maine\newline 72085 Le Mans Cedex,
France.}}

\maketitle

\begin{abstract}
A problem of goodness-of-fit test for ergodic diffusion processes is presented.  In the null hypothesis the drift of the diffusion is
supposed to be in a parametric form with unknown shift parameter. Two
Cramer-Von Mises type test statistics are studied. The first one is based on local time estimator of the
invariant density, the second one  is based on  the empirical distribution function. The
unknown parameter is estimated via the maximum likelihood estimator. It is shown
that both the limit distributions of the two test statistics do not depend on
the unknown parameter, so the distributions of the tests are asymptotically parameter free. Some considerations on the consistency of the proposed tests and some simulation studies are also given.

\end{abstract}

\pagestyle{fancy} \fancyhead{}
\fancyhead[C]{Negri,  Zhou}

\vspace{2mm}

\section{Introduction}
We consider the problem of goodness of fit test
for the model of ergodic diffusion process when
this process under the null hypothesis belongs to a given
parametric family. We study  the Cramer-von
Mises type statistics in two different cases. The first one is based on local time estimator and the second one is based on  empirical distribution
function estimator. We show that the  Cramer-von
Mises type statistics converge in both cases to some limits which do not depend on
the unknown parameter, so the test is asymptotically parameter
free (APF).

Let us remind the similar statement of the problem in the well known
case of the observations of independent identically distributed
random variables $X^n=\left(X_1,\ldots,X_n\right)$. Suppose that the
distribution of $X_j$ under hypothesis is $F\left(\vartheta
,x\right)=F\left({x-\vartheta}\right)$, where $\vartheta $ is some
unknown parameter. Then the Cramer-von Mises type test is
$$
\hat\psi_n \left(X^n\right)=\1_{\left\{\omega _n^2>e_\varepsilon
  \right\}},\qquad \omega _n^2=n\int_{-\infty }^{\infty }\left[\hat
  F_n\left(x\right)-F\left(x-\hat \vartheta _n\right)\right]^2{\rm
  d}F\left(x- \hat \vartheta _n\right)
$$
where the statistic $ \omega _n^2$ under hypothesis converges in
distribution to a random variable $\omega ^2$ which  does not depend
on $\vartheta $. Therefore the threshold $e_\varepsilon $ can
calculated as solution of the equation
$$
\P\left\{\omega ^2>e_\varepsilon \right\}=\varepsilon .
$$
The details concerning this result can be found in Darling \cite{D1955}.
For more general problems see the works of Kac, Kiefer \&
Wolfowitz \cite{KKW1955}, Durbin \cite{Durbin} or Martynov \cite{Martynov1979},
\cite{Martynov1992}.

A similar problem exists for the continuous time stochastic processes,
which are widely used as mathematic models in many fields. The
goodness of fit tests (GoF) are studied by many authors. For example
Kutoyants \cite{Kuto} discusses some possibilities of the
construction of such tests.  In particular,  he considers the
Kolmogorov-Smirnov statistics and the Cramer-von Mises Statistics
 based on the continuous observation.
Note that the Kolmogorov-Smirnov statistics for ergodic diffusion
process was studied in Fournie \cite{Fournie1992} and in Fournie and
Kutoyants \cite{FK1993}. However, due to the structure of the
covariance of the limit process, the Kolmogorov-Smirnov statistics
is not asymptotically distribution free in diffusion process
models. More recently Kutoyants  \cite{Kutoyants2010} has proposed
a modification of the Kolmogorov-Smirnov statistics for diffusion
models that became asymptotically distribution free. See also
Dachian and Kutoyants \cite{DK} where they propose some GoF tests
for diffusion and inhomogeneous Poisson processes with simple basic
hypothesis. It was shown that these tests are asymptotically
distribution  free. In the case of Ornstein-Uhlenbeck process
Kutoyants showed that the Cramer-von Mizes type tests are
asymptotically parameter free \cite{Kutoyants2012}. Another test was
studied by Negri and Nishiyama \cite{NN}.

\section{Main Results}

 Suppose that we observe an ergodic diffusion process, solution to the
 following stochastic differential equation
\begin{equation}
\label{equ} \d X_t=S(X_t)\d t+\d W_t,\quad X_0,\ 0\leq t\leq T.
\end{equation}

We want to test the following null hypothesis
$$
\mathcal{H}_0\quad :\qquad S\left(x\right)=S_*\left(x-\vartheta
\right),\quad \vartheta \in \Theta,
$$
where $S_*\left(\cdot \right)$ is some known function and the shift
parameter $\vartheta $  is unknown. We suppose that $0\in\Theta=\left(\alpha ,\beta \right)$.
Let us introduce the family
$$
{\mathcal S}\left(\Theta \right)=\left\{S_*\left(x-\vartheta
\right),\quad \vartheta \in \Theta=(\alpha,\beta)\right\}.
$$
The alternative is defined as
$$
\mathcal{H}_1\quad :\qquad S\left(\cdot\right)\not \in
\overline{\mathcal {S}(\Theta)},
$$
where  $\overline{\mathcal {S}(\Theta)}=\left\{S\left(x-\vartheta
\right),\vartheta \in \left[\alpha ,\beta \right] \right\} $.

We suppose that the trend coefficients $S\left(\cdot \right)$ of the
observed diffusion process under both hypotheses satisfy the
conditions:

$\mathcal{E}\mathcal{S}$. {\it The function $S(\cdot)$ is locally bounded
and for some $C>0$, }
$$
xS(x)\leq C(1+x^2).
$$
and\\

$\mathcal{A}_0$. {\it The function $S(\cdot)$ satisfies}
\begin{equation}
  \label{conditionA}
     \varlimsup\limits_{|x|\rightarrow\infty}{\rm
    sgn}(x)S(x)<0.
\end{equation}

Remind that under the condition $\mathcal{E}\mathcal{S}$, the
equation \eqref{equ} has a unique weak solution  (See
\cite{Durett}). Moreover under the condition $\mathcal{A}_0$, the
diffusion process is recurrent and its invariant density $f(x,\vth)$
under hypothesis $\mathcal{H}_0$ can be given explicitly (See
\cite{Kuto}, Theorem 1.16):
$$
f(x,\vth)=\frac{1}{G(\vth)}\exp\left\{2\int_\vth^xS_*(y-\vth)\d
y\right\}.
$$

Denote by $\xi_\vth$ a random variable (r.v.) having this density and
the corresponding mathematic expectation by $\E_\vth$. To simplify
the notations, for the case $\vth=0$, we denote  the density function as
$f(x)=f(x,0)$, and the corresponding distribution function as $F(x)$;
correspondingly the r.v. is $\xi_0$, and the mathematical expectation
is $\E_0$. Denote $\mathcal P$ as the class of functions having
polynomial majorants i.e.
\[\mathcal P=\{h(\cdot):\ |h(x)|\leq C(1+|x|^p)\},\]
with some $p>0$. Let $h'(x)$ the derivative of $h(x)$ w.r.t. $x$. \bigskip

Let us fix some $\vep\in(0,1)$, and  denote by $\mathcal {K}_\vep$
the class of tests $\psi_T$ of asymptotic size $\vep$, i.e.
$$
\E_0\psi_T=\vep+o(1).
$$
Our object is to construct this kind of tests.\\

To verify the hypothesis $\mathcal H_0$, we propose two tests. The
first one is based on the local time estimator (LTE) $\hat f_T(x)$
of the invariant density, which can be written as
$$
\hat f_T(x)=\frac{1}{T}(|X_T-x|-|X_0-x|)-\frac{1}{T}\int_0^T{\rm
sgn}(X_t-x)\d X_t.
$$
The unknown parameter is estimated via the maximum likelihood estimator (MLE)
$\hat\vth_T$, which is defined as the solution of the equation
$$
L(\hat\vth_T,X^T)=\sup\limits_{\theta\in\Theta}L(\theta,X^T),
$$
where $L(\vth,X^T)$ is the log-likelihood ratio
$$
L(\vth, X^T)=\int_0^TS_*(X_t-\vth)\d
X_t-\frac{1}{2}\int_0^TS_*(X_t-\vth)^2\d t. $$

We give the following regularity conditions $\mathcal{A}$ to have
the consistency and the asymptotical normality of the MLE:\\

\textbf{Condition $\mathcal {A}$.}

 {\it $\mathcal{A}_1$. The function $S_*(\cdot)$ is continuously
differentiable, the derivative $S'_*(\cdot)\in\mathcal P$ and is
uniformly continuous in the following sense:}

\[\lim\limits_{\nu\rightarrow0}\sup\limits_{|\tau|<\nu}
\E_0\big|S'_*(\xi_0)-S'_*(\xi_0+\tau)\big|^2=0.\]

{\it$\mathcal{A}_2$. The Fisher information}
\begin{equation}
\label{Fisher}
  I=\E_0S'_*(\xi_0)^2>0.
\end{equation}
 {\it Moreover, for any $\nu>0$}
\[\inf\limits_{|\tau|>\nu}\E_0\big(S_*(\xi_0)-S_*(\xi_0+\tau)\big)^2>0.\]

Denote the statistic based on the LTE as follows
$$
\delta_T=T\int_{-\infty}^\infty\left(\hat
f_T(x)-f(x-\hat\vth_T)\right)^2\d x,
$$
we will prove that under hypothesis $\mathcal {H}_0$, it converges
in distribution to
\begin{equation}
\label{limDelta}
\delta=\int_{-\infty}^\infty\left(\int_{-\infty}^\infty
\left(2f(x)\frac{\1_{\{y>x\}}-F(y)}{\sqrt{f(y)}}-\frac{1}{I}S_*'(y)\sqrt{f(y)}f'(x)\right)\d
W(y)\right)^2\d x,
\end{equation}
with $W(y)=W_1(y),\ y\in \mathbb{R}^+$, $W(y)=W_2(-y),\ y\in
\mathbb{R}^-$, where $W_1$ and $W_2$ are independent Wiener
processes. The Cramer-von Mises type test is defined as
$$
\psi_T=\1_{\{\delta_T>d_\vep\}},
$$
where $d_\vep$ is the $1-\vep$ quantile of the distribution of
$\delta$, that is the solution of the following equation
\begin{equation}
\label{d_vep}\P\Big(\delta\geq d_\vep\Big)=\vep.
\end{equation}

The main result for the Cramer von Mises test based on local time estimator is the following:
\begin{theorem}
\label{MainResult1}
  Let the conditions $\mathcal {E}\mathcal {S}$, $\mathcal {A}_0$ and $\mathcal
  {A}$ be fulfilled, then the test {\rm
  $\psi_T=\1_{\{\delta_T>d_\vep\}}$}
  belongs to $\mathcal {K}_\vep$.
\end{theorem}
The theorem is proved in  Section \ref{sectionLTE}.

Note that neither $\delta$ nor $d_\vep$ depends on the unknown
parameter. This allows us to  conclude that the test is APF.\\

The second test is based on the same MLE and the empirical
distribution function (EDF):
\[\hat F_T(x)=\frac{1}{T}\int_{0}^T\1_{\{X_t<x\}}\d t.\]
The corresponding statistic is
\[\Delta_T=T\int_{-\infty}^\infty\left(\hat F_T(x)-F(x-\hat\vth_T)\right)^2\d x,\]
which converges in distribution to
\begin{equation}
\label{limdelta}
\Delta=\int_{-\infty}^\infty\left(\int_{-\infty}^\infty
\left(2\frac{F(y\wedge
x)-F(y)F(x)}{\sqrt{f(y)}}-\frac{1}{I}S_*'(y)\sqrt{f(y)}f(x)\right)\d
W(y)\right)^2\d x.
\end{equation}
Thus we propose the Cramer-von Mises type test
$$
\Psi_T=\1_{\{\Delta_T>c_\vep\}},
$$
where $c_\vep$ is the solution of the equation
\begin{equation}
\label{c_vep}\P\Big(\Delta\geq c_\vep\Big)=\vep.
\end{equation}

The main result for the Cramer von Mises test based on empirical distribution function estimator is the following:
\begin{theorem}
\label{MainResult2}
  Under conditions $\mathcal {E}\mathcal {S}$, $\mathcal {A}_0$ and $\mathcal
  {A}$, the test {\rm$\Psi_T=\1_{\{\Delta_T>c_\vep\}}$}
  belongs to $\mathcal {K}_\vep$.
\end{theorem}
The theorem is proved In Section \ref{sectionEDF}.

\bigskip

\section{Proof of Theorem \ref{MainResult1}}
\label{sectionLTE}

In this section, we study the test
$\psi_T=\1_{\{\delta_T>d_\vep\}}$,  where
$$\delta_T=T\int_{-\infty}^\infty\left(\hat f_T(x)-f(x-\hat\vth_T)\right)^2\d x.$$

Under the basic hypothesis $\mathcal {H}_0$, the density of the
invariant law can be presented as follows:
\begin{eqnarray*}
  f(x,\vth)
  &=&\frac{\exp\{2\int_{\vth}^xS_*(y-\vth)\d y\}}{\int_{-\infty}^\infty\exp\{2\int_{\vth}^yS_*(z-\vth)\d z\}\d y}\\
  &=&\frac{\exp\{2\int_{0}^{x-\vth}S_*(y)\d y\}}{\int_{-\infty}^\infty\exp\{2\int_{0}^{y-\vth}S_*(z)\d z\}\d y}\\
  &=&f(x-\vth).
\end{eqnarray*}
Note that the distribution function of the process satisfies
\begin{eqnarray*}
  F(x,\vth)
  &=&\int_{-\infty}^x f(y-\vth)\d y
  =\int_{-\infty}^{x-\vth} f(y)\d y=F(x-\vth).
\end{eqnarray*}
In addition, for any integrable function $h$,
\begin{eqnarray}
\label{TransformE}
   \nonumber
   \E_{\vth}h(\xi_\vth-\vth)&=&\int_{-\infty}^\infty
   h(x-\vth)f(x-\vth)\d x\\
   &=&\int_{-\infty}^\infty h(x)f(x)\d x=\E_0
   h(\xi_0).
\end{eqnarray}
Note that the Fisher information in our case does not depend on the
unknown parameter $\vth$:
$$
I=\E_{\vth_0}S_*'(\xi_{\vth_0}-\vth_0)^2=\E_0 S_*'(\xi_0)^2>0.
$$
where $\vth_0$ is the true value of the unknown parameter.

From the condition $\mathcal {A}_0$, it follows that there exist
some constants $A>0$ and $\gamma>0$ such that for all $|x|>A$,
\begin{equation}
  \label{extremeS}
  {\rm sgn}(x)S_*(x)<-\gamma.
\end{equation}
It can be shown that
 for $x>A$,
\begin{eqnarray*}
  f(x)&=&\frac{1}{G(S_*)}{\rm exp}\left\{2\left(\int_0^A
  +\int_A^{x}\right)S_*(y)\d y\right\}
  <C{\rm e}^{-2\gamma x}.
\end{eqnarray*}\\
Similar result can be deduced for $x<-A$,
so we have
\begin{equation}\label{extremef}
  f(x)<C{\rm e}^{-2\gamma |x|},\quad{\rm for}\ |x|>A.
\end{equation}

Let the conditions $\mathcal{A}_0$ and $\mathcal{A}$ be fulfilled,
  then the MLE $\hat\vth_T$ is consistent, i.e., for any $\nu>0$,
  {\rm\[\lim\limits_{T\rightarrow\infty}\P_{\vth_0}\big\{|\hat\vth_T-\vth_0|>\nu\big\}=0;\]}
 it is  asymptotically normal
\begin{equation}
  \mathcal{L}_{\vth_0}\big\{\sqrt T(\hat\vth_T-\vth_0)\big\}\Longrightarrow\mathcal
  N(0,I^{-1});
  \label{lemmaK}
  \end{equation}
  and the moments converge i.e., for $p>0$
  {\rm$$
  \lim\limits_{T\rightarrow\infty}\E_{\vth_0}\left|\sqrt T(\hat \vth_T-\vth_0)\right|^p
  =\E_0\left|\hat u\right|^p,
  $$}
where $\hat u\sim\mathcal
N(0,I^{-1})$.
The proof can be found in \cite{Kuto},Theorem 2.8.
We can define
\[
\hat u=\frac{1}{I}\int_{-\infty}^\infty S_*'(y)\sqrt{f(y)}\d W(y),
\]
and denoted $\hat u_T=\sqrt{T}(\hat\vth_T-\vth_0)$, the asymptotical
normality \eqref{lemmaK} can be written as
\begin{equation}
\label{convMLE}
  \mathcal{L}_{\vth_0}\left\{\hat u_T\right\}\Longrightarrow\mathcal{L}\left\{\hat
u\right\}.
\end{equation}
\\

We define $\eta_T(x)=\sqrt T\left(\hat f_T(x)-f(x-\vth_0)\right)$. In
\cite{Kuto} Theorem 4.11, we can find  the following representation
\begin{eqnarray}
  \nonumber
  \eta_T(x)&=&\sqrt T(\hat f_T(x)-f(x-\vth_0))\\
  \nonumber
  &=&2\frac{f(x-\vth_0)}{\sqrt T}\int_{X_0}^{X_T}\left(\frac{\1_{\{y>x\}}-F(y-\vth_0)}{f(y-\vth_0)}\right)\d
  y\\
  &&-2\frac{f(x-\vth_0)}{\sqrt T}\int_0^T\left(\frac{\1_{\{X_t>x\}}-F(X_t-\vth_0)}{f(X_t-\vth_0)}\right)\d
  W_t.\label{lemmaLTEEqu1}
\end{eqnarray}
Let us put
$$
M(y,x)=2f(x)\frac{\1_{\{y>x\}}-F(y)}{f(y)}.
$$
Then $\eta_T(x)$ can be written as
\begin{eqnarray}
  \nonumber
  \eta_T(x)&=&
  \frac{1}{\sqrt T}\int_{X_0}^{X_T}M(y-\vth_0,x-\vth_0)\d
  y\\
  &&-\frac{1}{\sqrt T}\int_0^TM(X_t-\vth_0,x-\vth_0)\d
  W_t.\label{lemmaLTEEqu1-1}
\end{eqnarray}
We can state
\begin{lemma}\label{IntLim}Let the condition $\mathcal {A}_0$ be
fulfilled, then
  {\rm$$
  \int_{-\infty}^\infty\E_0\left(\int_0^{\xi_0}M(y,x)\d y\right)^2\d
  x<\infty.
  $$}
\end{lemma}
\noindent\textbf{Proof.} Applying the estimate \eqref{extremef}, for
$x>A$,
\begin{eqnarray*}
  &&\E_0\left(\int_{0}^{\xi_0}M(y,x)\d y\right)^2\\
  &&\quad=4f(x)^2\int_{-\infty}^{\infty}\left(\int_0^z\frac{\1_{\{y>x\}}-F(y)}{f(y)}\d
  y\right)^2f(z)\d z\\
  &&\quad=4f(x)^2\left(\int_{-\infty}^{-A}+\int_{-A}^{A}+\int_{A}^{x}
  \right)\left(\int_0^z\frac{-F(y)}{f(y)}\d
  y\right)^2f(z)\d z\\
  &&\quad\quad+4f(x)^2\int_{x}^{\infty}\left(\int_0^{x}\frac{-F(y)}{f(y)}\d
  y+\int_{x}^z\frac{1-F(y)}{f(y)}\d y\right)^2f(z)\d z
\end{eqnarray*}
Further,
\begin{eqnarray*}
  &&f(x)^2\int_{-\infty}^{-A}\left(\int_0^z\frac{-F(y)}{f(y)}\d
  y\right)^2f(z)\d z\\
  &&\quad=f(x)^2\int_{-\infty}^{-A}\left(\left(\int_z^{-A}
  +\int_{-A}^0\right)\frac{F(y)}{f(y)}\d y\right)^2f(z)\d z\\
  &&\quad\leq f(x)^2\int_{-\infty}^{-A}\left(\int_z^{-A}\int_{-\infty}^y
  \frac{1}{G}{\rm exp}\left(-2\int_{u}^yS_*(v)\d v\right)\d u\d y+C_1\right)^2f(z)\d z\\
  &&\quad\leq f(x)^2\int_{-\infty}^{-A}\left(C_2\int_z^{-A}\int_{-\infty}^y
  {\rm e}^{-2\gamma(y-u)}\d u\d y+C_1\right)^2f(z)\d z\\
  &&\quad\leq C f(x)^2\int_{-\infty}^{-A}(1+z)^2f(z)\d z\leq C f(x)^2\leq C{\rm e}^{-4\gamma x},
\end{eqnarray*}
moreover
\begin{eqnarray*}
  &&f(x)^2\int_A^x\left(\int_0^{z}\frac{-F(y)}{f(y)}\d y\right)^2f(z)\d
  z\\
  &&\quad \leq\int_A^x\left(\left(\int_0^{A}+\int_A^z\right)\frac{f(x)}{f(y)}\d y\right)^2f(z)\d
  z\\
  &&\quad\leq\int_A^x\left(C_1f(x)+C_2\int_A^z{\rm e}^{-2\gamma (x-y)}\d y\right)^2f(z)\d
  z\\
  &&\quad\leq \int_A^x\left(C_1{\rm e}^{-2\gamma x}+C_2'{\rm e}^{-2\gamma (x-z)}
  -C_2'{\rm e}^{-2\gamma (x-A)}\right)^2\cdot C{\rm e}^{-2\gamma z}\d
  z\\
  &&\quad\leq{\rm e}^{-4\gamma x}\int_A^x
  \left(C_3{\rm e}^{2\gamma z}+C_4{\rm e}^{-2\gamma z}\right)\d z\leq C{\rm e}^{-2\gamma x},
\end{eqnarray*}
and finally
\begin{eqnarray*}
  &&f(x)^2\int_{x}^{\infty}\left(\int_x^{z}\frac{1-F(y)}{f(y)}\d y\right)^2f(z)\d
  z\\
  &&\quad\leq C f(x)^2\int_{x}^{\infty}\left(\int_x^{z}
  \int_y^\infty{\rm e}^{-2\gamma (u-y)}\d u\d y\right)^2{\rm e}^{-2\gamma z}\d
  z\\
  &&\quad\leq C f(x)^2\int_{x}^{\infty}(z-x)^2{\rm e}^{-2\gamma z}\d
  z\\
  &&\quad\leq C f(x)^2\int_{0}^{\infty}s^2{\rm e}^{-2\gamma (s+x)}\d
  s\leq C {\rm e}^{-6\gamma x}.
\end{eqnarray*}
Then we have
\begin{equation}
\label{extremeM-2}
  \E_0\left(\int_{0}^{\xi_0}M(y,x)\d y\right)^2
  \leq C{\rm e}^{-2\gamma |x|}\quad{\rm for}\ x>A.
\end{equation}
Similar estimate  can be obtained for $x<-A$, therefore the result holds
for $|x|>A$. We obtain finally
\begin{eqnarray*}
  &&\int_{-\infty}^\infty\E_0\left(\int_0^{\xi_0}M(y,x)\d y\right)^2\d
  x\\
  &&\quad=\left(\int_{-\infty}^{-A}+\int_{-A}^A+\int_A^\infty\right)\E_0\left(\int_0^{\xi_0}M(y,x)\d y\right)^2\d
  x\\
  &&\quad\leq C_1\int_{-\infty}^{-A}{\rm e}^{2\gamma x}\d x+C_2+C_3\int_{A}^{\infty}{\rm e}^{-2\gamma
  x}\d x
  <\infty.
\end{eqnarray*}
\bigskip

This result yields directly the conditions $\mathcal {O}$
of Theorem 4.11 in \cite{Kuto}:
$$
  \E_{\vth_0}M(\xi_{\vth_0}-\vth_0,x-\vth_0)^2=\E_0
  M(\xi_0,x-\vth_0)^2<\infty,
$$
and
$$
 \E_{\vth_0}\left(\int_{0}^{\xi_{\vth_0}}M(y-\vth_0,x-\vth_0)\d
 y\right)^2<\infty.
$$
So we can deduce  the convergence and the asymptotical normality of
$\eta_T(x)$. In fact  under the condition $\mathcal{A}_0$, the LTE $\hat
f_T(x)$ is  consistent and asymptotically normal, that is
{\rm\[\eta_T(x)=\sqrt T\left(\hat
f_T(x)-f(x-\vth_0)\right)\Longrightarrow\eta(x-\vth_0) ,\]}
 where
$\eta(x)\sim\mathcal {N}(0,d(x)^2)$, and
{\rm\[d(x)^2=4f(x)^2\E_0\left(\frac{\1_{\{\xi_0>x\}}-F(\xi_0)}{f(\xi_0)}\right)^2.\]}
Moreover
\begin{eqnarray*}
&&\E_{\vth_0}\left(\eta_T(x)\eta_T(y)\right)\\
&&\quad=4f(x-\vth_0)f(y-\vth_0)\E_0\left(\frac{\left(\1_{\{\xi_0>x-\vth_0\}}-F(\xi_0)\right)
\left(\1_{\{\xi_0>y-\vth_0\}}-F(\xi_0)\right)}{f(\xi_0)^2}\right).
\end{eqnarray*}

We can define
\[\eta(x)=\int_{-\infty}^\infty
M(y,x)\sqrt {f(y)}\d W(y).
\]
The distribution of $\eta(x)$ is $\mathcal N(0,\E_0 M(\xi_0,x)^2)$,
and we have the following convergence
\begin{equation}
\label{convLTE}
  \eta_T(x)\Longrightarrow\eta(x-\vth_0).\\
\end{equation}

For $\hat u_T$ and $\eta_T(x)$, we need more than \eqref{convMLE}
and convergence \eqref{convLTE}.
\begin{lemma}
\label{LemConv}
  Let conditions  $\mathcal {A}_0$ and  $\mathcal {A}$ be fulfilled,
  then $(\eta_T(x_1),...,\eta_T(x_k),\hat u_T)$ is asymptotically normal:
  {\rm\[\mathcal{L}\left(\eta_T(x_1),...,\eta_T(x_k),\hat u_T\right)\Longrightarrow
  \mathcal{L}\left(\eta(x_1-\vth_0),...,\eta(x_k-\vth_0),\hat u\right),\]}
  for any $\mathbf{x}=\{x_1,x_2,...,x_{k}\}\in\mathbb{R}^k$.
\end{lemma}
\noindent\textbf{Proof.} The first integral in
\eqref{lemmaLTEEqu1-1} converges to zero, so we only need to verify the convergence for
the part of It\^o integral. Let us denote for simplicity
\begin{eqnarray*}
  &&\eta_T^0(x)=
  \frac{1}{\sqrt T}\int_0^TM(X_t-\vth_0,x)\d
  W_t.
\end{eqnarray*}
It is sufficient to verify that for any $\mathbf{x}=\{x_1,x_2,...,x_{k}\}$,
\begin{equation}
  \label{LemConvEqu1}
  \left(\eta_T^0(x_1),...,\eta_T^0(x_k),\hat u_T\right)\Longrightarrow
  \left(\eta(x_1),...,\eta(x_k),\hat u\right).
\end{equation}
Remember that $\hat u_T$ can be defined as follows,
\begin{equation}
  \label{LemConvEqu2}
  Z_T(\hat u_T)=\sup\limits_{u\in\mathbbm{U}_T} Z_T(u),\quad
  \mathbbm{U}_T=\{u:\vth+\frac{u}{\sqrt T}\in\Theta\},
\end{equation}
where
\[
  Z_T(u)=\frac{\d \P_{\vth+\frac{u}{\sqrt T}}^T}{\d \P_\vth^T}(X^T)=
  {\rm exp}\left\{u\Lambda_T-\frac{u^2}{2}I+r_T\right\}.
\]
Here $\Lambda_T=\frac{1}{\sqrt T}\int_0^TS'_*(X_t-\vth_0)\d W_t$ and
$r_T\longrightarrow0$. It was proved in \cite{Kuto}, Theorem 2.8
that $Z_T(\cdot)$ converges in distribution to $Z(\cdot)$, where
$$Z(u)={\rm exp}\left\{u\Lambda-\frac{u^2}{2}I\right\},$$
where $\Lambda$ is a r.v. with normal distribution $\mathcal {N}(0,I)$,
which can be written as
\[\Lambda=\int_{-\infty}^\infty S_*'(y)\sqrt{f(y)}\d W(y).\]
Therefore
\[\hat u_T\Longrightarrow \hat u=\frac{\Lambda}{I}.\]
Take $\mathbf{u}=\{u_1,u_2,...,u_{m}\}$. We have to verify that the
joint finite-dimensional distribution of $Y_T$

\[Y_T=\left(\eta_T^0(x_1),\eta_T^0(x_2),...,\eta_T^0(x_k),Z_T(u_1),Z_T(u_2),...,Z_T(u_m)\right)\]
converges to the finite-dimensional distribution of $Y$
\[Y=\left(\eta(x_1),\eta(x_2),...,\eta(x_k),Z(u_1),Z(u_2),...,Z(u_m)\right).\]
Note that the only stochastic term in $Z_T(u)$ is $\Lambda_T$, so
\eqref{LemConvEqu1} is equivalent  to
\begin{equation}
  \label{LemConvEqu4}
  \left(\eta_T^0(x_1),\eta_T^0(x_2),...,\eta_T^0(x_k),\Lambda_T\right)\Longrightarrow
  \left(\eta(x_1),\eta(x_2),...,\eta(x_k),\Lambda\right).
\end{equation}
Take $\mathbf{\lambda}=\{\lambda_1,\lambda_2,...,\lambda_{k+1}\}$,
and put
$$h(y,\mathbf{x},\mathbf{\lambda})=\sum\limits_{l=1}^{k}\lambda_lM(y,x_l)
  +\lambda_{k+1} S_*'(y).$$
We have
\begin{eqnarray*}
  &&\E_{\vth_0}h(\xi_{\vth_0}-\vth_0,\mathbf{x},\mathbf{\lambda})^2
  =\E_0 h(\xi_0,\mathbf{x},\mathbf{\lambda})^2\\
  &&\quad=\int_{-\infty}^\infty\left(\sum\limits_{l=1}^{k}\lambda_lM(y,x_l)
  +\lambda_{k+1} S_*'(y)\right)^2f(y)\d y\\
  &&\quad=\int_{-\infty}^\infty\left(\sum\limits_{l=1}^{k}
  2\lambda_lf(x_l)\frac{\1_{\{y>x_l\}}-F(y)}{\sqrt{f(y)}}
  +\lambda_{k+1} S_*'(y)\sqrt{f(y)}\right)^2f(y)\d y\\
  &&\quad=\int_{-\infty}^\infty\left(\sum\limits_{l=1}^{k}
  \sum\limits_{m=1}^{k}4\lambda_l\lambda_mf(x_l)f(x_m)
  \frac{(\1_{\{y>x_l\}}-F(y))(\1_{\{y>x_m\}}-F(y))}{f(y)}\right.\\
  &&\quad\quad\left.+\sum\limits_{l=1}^{k}\lambda_l\lambda_{k+1}\left(\1_{\{y>x_l\}}-F(y)\right)S_*'(y)
  +\lambda_{k+1}^2S_*'(y)^2f(y)\right)\d y
   <\infty.
\end{eqnarray*}
The law of large number gives us
\[
  \frac{1}{T}\int_0^Th(X_t-\vth_0,\mathbf{x},\mathbf{\lambda})^2\d t
  \longrightarrow\E_0 h(\xi_0,\mathbf{x},\mathbf{\lambda})^2.
\]
Moreover, the central limit theorem for stochastic integral gives us

$$\frac{1}{\sqrt T}\int_0^Th(X_t-\vth_0,\mathbf{x},\mathbf{\lambda})\d W_t
  \Longrightarrow\mathcal {N}\left(0,\E_0 h(\xi_0,\mathbf{x},\mathbf{\lambda})^2\right).$$
In addition
$\sum\limits_{l=1}^k\lambda_l\eta(x_l)+\lambda_{k+1}\Lambda$ is a
zero mean normal r.v. with variance
\begin{eqnarray*}
  &&\E_0 \left(\sum\limits_{l=1}^k\lambda_l\eta(x_l)+\lambda_{k+1}\Lambda\right)^2\\
  &&\quad=\sum\limits_{l=1}^k\sum\limits_{m=1}^k\lambda_l\lambda_m\E_0\left(\eta(x_l)\eta(x_m)\right)
  +\sum\limits_{l=1}^k\lambda_l\lambda_{k+1}\E_0(\eta(x_l)\Lambda)+\lambda_{k+1}^2\E_0(\Lambda)^2.\\
\end{eqnarray*}
Furthermore
\begin{eqnarray*}
  &&\E_0\left(\eta(x_l)\eta(x_m)\right)\\
  &&\quad=4f(x_l)f(x_l)
  \int_{-\infty}^\infty\frac{(\1_{\{y>x_l\}}-F(y))(\1_{\{y>x_m\}}-F(y))}{f(y)}\d
  y,
\end{eqnarray*}
and
\begin{eqnarray*}
  \E_0(\eta(x_l)\Lambda)=-2f(x_l)\int_{-\infty}^\infty(\1_{\{y>x_l\}}-F(y))S_*'(y)\d y,
\end{eqnarray*}
\begin{eqnarray*}
  \E_0(\Lambda)^2=\int_{-\infty}^\infty S_*'(y)^2f(y)\d y.
\end{eqnarray*}
We find that
\[\E_{\vth_0} h(\xi_{\vth_0}-\vth_0,\mathbf{x},\mathbf{\lambda})^2
=\E_0 h(\xi_0,\mathbf{x},\mathbf{\lambda})^2=
\E_0\left(\sum\limits_{l=1}^k\lambda_l\eta(x_l)+\lambda_{k+1}\Lambda\right)^2.\]
 This is as to say
\[
\sum\limits_{l=1}^k\lambda_l\eta_T^0(x_l)+\lambda_{k+1}\Lambda_T\Longrightarrow
\sum\limits_{l=1}^k\lambda_l\eta(x_l)+\lambda_{k+1}\Lambda
\]
thus  \eqref{LemConvEqu1} follows from this last convergence in distribution, and so
the lemma is proved.\\
\bigskip

\begin{lemma}
\label{LemIntConv}
  Let conditions $\mathcal {A}_0$ and $\mathcal {A}$ be fulfilled, then
  {\rm\[\mathcal{L}\left\{\int_{-\infty}^\infty\left(\eta_T^0(x)-\hat u_Tf'(x)\right)^2\d x\right\}\Longrightarrow
  \mathcal{L}\left\{\int_{-\infty}^\infty\left(\eta(x)-\hat u f'(x)\right)^2\d
  x\right\}\]}
\end{lemma}

\noindent\textbf{Proof.} Denote $\zeta_T(x)=\eta_T^0(x)-\hat
u_Tf'(x)$ and $\zeta(x)=\eta(x)-\hat u f'(x)$, we will prove the
following properties

i) For $x,y\in[-L,L]$ and $|x-y|\leq 1$,
\begin{equation}
\label{LemIntConvEqu2}
  \E_{\vth_0}|\zeta_T(x)^2-\zeta_T(y)^2|^2\leq C|x-y|^\delta,\quad {\rm\ with\ some\ }\delta>0.
\end{equation}

ii) $\forall \vep>0,\ \exists L>0$, such that
\begin{equation}
\label{LemIntConvEqu3}
  \E_{\vth_0}\int_{\{|x|>L\}}\zeta_T(x)^2\d
  x<\vep,\quad \forall T>0.
\end{equation}
From i) it follows  the convergence in every  bounded set $[-L, L]$:
  \[\mathcal{L}\big\{\int_{-L}^L\zeta_T(x)^2\d x\big\}\Longrightarrow
  \mathcal{L}\big\{\int_{-L}^L\zeta(x)^2\d x\big\}.\]
The result in i) along with ii) gives us the result. \\

First we prove i). We have
$$
  \E_{\vth_0}\left(\zeta_T(x)^2\right)\leq 2\E_{\vth_0}\eta_T^0(x)^2+2f(x)^2\E_{\vth_0}\hat u_T^2\leq C.
$$
\begin{eqnarray*}
  &&\E_{\vth_0}\left|\zeta_T(x)^2-\zeta_T(y)^2\right|^2\\
  &&\quad=\E_{\vth_0}\left(|\zeta_T(x)+\zeta_T(y)|^2|\zeta_T(x)-\zeta_T(y)|^2\right)\\
  &&\quad\leq C\E_{\vth_0}|\zeta_T(x)-\zeta_T(y)|^2\\
  &&\quad\leq C(f'(x)-f'(y))^2\E_{\vth_0}|\hat u_T|^2
  +\E_{\vth_0}|(\eta_T^0(x)-\eta_T^0(y))|^2.\\
\end{eqnarray*}
For the first part, let us recall the following result, given  in \cite{Kuto},
page 119: for any $p>0,\ R>0$, chosen $N$ sufficiently large, we have
\[\P_{\vth_0}^T\left\{|\hat u_T|^p>R\right\}\leq\frac{C_N}{R^{N/p}}.\]
Now, denoted $F_T(u)$ the distribution of $|\hat u_T|$, we have
\begin{eqnarray}
  \nonumber
  &&\E_{\vth_0}|\hat u_T|^p=\int_0^\infty u^p\d F_T(u)
  \leq1-\int_1^\infty u^p\d[1-F_T(u)]\\
  \label{u_Tbounded}
  &&\quad\leq1-[1-F_T(1)]+p\int_1^\infty u^{p-1}\frac{C_N}{u^{N/p}}\d u\leq
  C.
\end{eqnarray}
Remember that under condition $\mathcal {A}_1$, $S_*$ and $f$ are
sufficiently smooth. So, for $x,y \in[-L,L]$ we can write
\begin{eqnarray*}
|f(x)-f(y)|&=&|f'(z)(x-y)|=|2S_*(z)f(z)(x-y)|\leq C|x-y|,
\end{eqnarray*}
and
\[|f'(x)-f'(y)|=|f''(z)(x-y)|=\left|4f(z)S_*^2(z)+2f(z)S_*'(z)\right||x-y|\leq C|x-y|.\]
So we have
$$
(f'(x)-f'(y))^2\E_{\vth_0}|\hat u_T|^2\leq
C|x-y|^2.
$$
For the second part, we can write
\begin{eqnarray*}
  &&\E_{\vth_0}|(\eta_T^0(x)-\eta_T^0(y))|^2\\
  &&\quad= C_1\E_{\vth_0}\left(\frac{1}{\sqrt T}\int_0^{T}
  (M(X_t-\vth_0,x)-M(X_t-\vth_0,y))\d W_t\right)^2\\
  &&\quad\leq\frac{C_1}{T}\int_0^{T}
  \E_{\vth_0}\left(M(X_t-\vth_0,x)-M(X_t-\vth_0,y)\right)^2\d
  t\\
  &&\quad=C_1\E_0\left(M(\xi_0,x)-M(\xi_0,y)\right)^2.
\end{eqnarray*}
Suppose that $x\leq y$,
\begin{eqnarray*}
  &&\E_0\left(M(\xi_0,x)-M(\xi_0,y)\right)^2\\
  &&\quad=\int_{-\infty}^x\left(2\frac{F(z)}{f(z)}(f(x)-f(y))\right)^2f(z)\d z\\
  &&\qquad+\int_x^y\left(2\frac{1}{f(z)}\left((1-F(z))f(x)
   +F(z)f(y)\right)\right)^4f(z)\d z\\
  &&\qquad+\int_y^\infty \left(2\frac{1-F(z)}{f(z)}(f(x)-f(y))\right)^2f(z)\d z\\
  &&\quad\leq C_1(x-y)^4+C_2(x-y)+C_3(x-y)^2\leq C(y-x).
\end{eqnarray*}

Similar result holds for  $x > y$. Then we obtain
\[\E_{\vth_0}\left|\eta_T^0(x)^2-\eta_T^0(y)^2\right|^2\leq C|x-y|,\quad x,y\in\mathbb{R}.\]
Thus we have
\[\E_{\vth_0}\left|\zeta_T(x)^2-\zeta_T(y)^2\right|^2\leq C|x-y|.\]

Now we prove ii). As in Lemma \ref{IntLim}, we can deduce that
$$
\E_0M(\xi_0,x)^2\leq C{\rm e}^{-2\gamma x},\quad \textrm{for } x>A.
$$
So for $L>A$,
\begin{eqnarray*}
  &&\E_{\vth_0}\int_L^\infty\left(\eta_T^0(x)\right)^2\d
  x=\E_{\vth_0}\int_L^\infty\left(\frac{1}{\sqrt T}\int_0^TM(X_t-\vth_0,x)\d W_t\right)^2\d
  x\\
  &&\quad\leq C\int_L^\infty\E_0 M(\xi_0,x)^2\d x
  \leq C\int_L^\infty{\rm e}^{-2\gamma x}\d x\leq C{\rm e}^{-2\gamma
  L}.
\end{eqnarray*}
Note that $f'(x)=2S_*(x)f(x)$ and  along with \eqref{u_Tbounded} we get
\begin{eqnarray*}
  &&\int_L^\infty\E_{\vth_0}\big(\eta_T^0(x)-f'(x)\hat u_T\big)^2\d x\\
  &&\quad\leq\int_L^\infty\left(2\E_{\vth_0}\eta_T(x)^2+2f'(x)\E_{\vth_0}\hat u_T^2\right)\d x\\
  &&\quad\leq\int_L^\infty C{\rm e}^{-2\gamma x}\d x=C{\rm e}^{-2\gamma
  L}.
\end{eqnarray*}
For any $\vep>0$, take $L=-\frac{\ln(\vep/C)}{2\gamma}\vee A$, then  we
have \eqref{LemIntConvEqu3}.\\
\bigskip\\
\textbf{Proof of Theorem \ref{MainResult1}.}\\

We can write
\begin{eqnarray*}
  \delta_T&=&T\int_{-\infty}^\infty(\hat f_T(x)-f(x-\hat\vth_T))^2\d x\\
  &=&T\int_{-\infty}^\infty\left((\hat f_T(x)-f(x-\vth_0))+(f(x-\vth_0)-f(x-\hat\vth_T))\right)^2\d x\\
  &=&\int_{-\infty}^\infty\left(\sqrt{T}(\hat f_T(x)-f(x-\vth_0))-\sqrt{T}(\hat\vth_T-\vth_0)f'(x-\tilde\vth_T)\right)^2\d
  x\\
  &=&\int_{-\infty}^\infty\left(\eta_T(x)-\hat u_T f'(x-\tilde\vth_T)\right)^2\d
  x.\\
\end{eqnarray*}
See that
\begin{eqnarray*}
  &&\E_{\vth_0} \int_{-\infty}^\infty\left(\hat u_T^2|
  f'(x-\tilde\vth_T)-f'(x-\vth_0)|^2\right)\d x\\
  &&\quad=\E_{\vth_0} \int_{-\infty}^\infty\left(\hat u_T^2
  f''(x-\vth_T^*)^2(\tilde \vth_T-\vth_0)^2\right)\d x,\\
\end{eqnarray*}
and that $f'(x-\vth)=S_*(x-\vth)f(x-\vth)$,
$f''(x,\vth)=S_*'(x-\vth)f(x-\vth)+S_*(x-\vth)^2f(x-\vth)$, the
smoothness of $S_*(\cdot)$ gives us the convergence
\[\E_{\vth_0} \int_{-\infty}^\infty\left(\hat u_T^2|
  f'(x-\tilde\vth_T)-f'(x-\vth_0)|^2\right)\d x\longrightarrow0.\]
\vspace{3mm}\\
Applying Lemma \ref{IntLim} and Lemma \ref{LemIntConv} we get
 \begin{eqnarray*}
  \delta_T&=&\int_{-\infty}^\infty\left(\eta_T^0(x-\vth_0)
  -\hat u_T f'(x-\vth_0)\right)^2\d x+o(1)\\
  &\Longrightarrow&\int_{-\infty}^\infty\left(\eta(x-\vth_0)
  -\hat u f'(x-\vth_0)\right)^2\d x\\
  &=&\int_{-\infty}^\infty\left(\eta(y)
  -\hat u f'(y)\right)^2\d y=\delta.
\end{eqnarray*}
\bigskip
We see that the limit of the statistic $\delta$ does not depend on
$\vth_0$, and the test $\psi_T=\1_{\{\delta_T\geq d_\vep\}}$ with
$d_\vep$ defined by
\[\P\Big(\delta\geq d_\vep\Big)=\vep\]
belongs to $\mathcal {K}_\vep$.\\

The same  procedure  can be applied with other estimators of the
  unknown parameter and of the invariant density, provided that they are consistent and asymptotically
  normal.
  For example, we can  take the minimum distance estimator (MDE)
  $\vth_T^*$ for $\vth_0$:
  $$
  \vth_T^*={\rm arg}\inf\limits_{\theta\in\Theta}\|\hat F(\cdot)-F(\theta,\cdot)\|,
  $$
  and the kernel estimators $\bar f_T(x)$ as estimator for the invariant density
 \[\bar f_T(x)=\frac{1}{\sqrt T}\int_0^TK(\sqrt T(X_t-x))\d t.\]
 Under some regularity conditions, the MDE
 $\hat\vth_T^*$ is asymptotically normal (See \cite{FK1993} or \cite{Kuto}):
 $$
 u^*_T=\sqrt T(\vth_T^*-\vth_0)\Longrightarrow \hat
u^*\sim\mathcal {N}(0,R(\vth_0)).
$$
Also if we do not present explicitly $R(\cdot)$ here, it  can be verified that $R(\vth)=R(0)$ does not depend on $\vth$. The kernel estimator
$\bar f_T(x)$ has the same asymptotic properties of the LTE (See
\cite{Kuto}). Then we can construct the statistic
$$
\mu_T=T\int_{-\infty}^\infty\left(\bar f(x)-f(x-
\vth_T^*)\right)^2\d x,
$$
which converges to
$$
\mu=\int_{-\infty}^\infty\left(\eta(x)- u^*f'(x)\right)^2\d x,
$$
that does not depend on the unknown parameter. So that the test
${\rm \1}_{\{\mu_T>k_\vep\}} $ with $k_\vep$ the solution of the
equation
$$
{\rm\P}\left(\mu>k_\vep\right)=\vep
$$
belongs to $\mathcal {K}_\vep$.

\section{Proof of Theorem \ref{MainResult2}}
\label{sectionEDF}

In this section, we study the GoF test $\Psi_T=\1_{\{\Delta_T\geq
c_\vep\}}$ defined by the  statistic
\[\Delta_T=T\int_{-\infty}^\infty\big(\hat F_T(x)-F(x-\hat\vth_T)\big)^2\d x,\]
where $\hat F_T(x)$ is the empirical distribution function:
\[\hat F_T(x)=\frac{1}{T}\int_0^T\1_{\{X_t<x\}}\d t.\]

Denote $\eta_T^F(x)=\sqrt{T}(\hat F_T(x)-F(x-\vth_0))$ and
$$
  H(z,x)=2\frac{F(z\wedge x)-F(z)F(x)}{f(z)}.
$$
In \cite{Kuto} Theorem 4.6,  the following equality is presented:
\begin{eqnarray*}
  \eta_T^F(x)
  &=&\frac{2}{\sqrt T}\int_{X_0}^{X_T}\frac{F((z\wedge x)-\vth_0)
  -F(z-\vth_0)F(x-\vth_0)}{f(z-\vth_0)}\d z\\
  &&\quad-\frac{2}{\sqrt{T}}\int_{0}^{T}\frac{F((X_t\wedge x)-\vth_0)
  -F(X_t-\vth_0)F(x-\vth_0)}{f(X_t-\vth_0)}\d W_t.
\end{eqnarray*}
Then
\begin{eqnarray*}
  \eta_T^F(x)&=&\frac{2}{\sqrt T}\int_{X_0}^{X_T}\frac{F((z-\vth_0)\wedge (x-\vth_0))
  -F(z-\vth_0)F(x-\vth_0)}{f(z-\vth_0)}\d z\\
  &&-\frac{2}{\sqrt{T}}\int_{0}^{T}\frac{F((X_t-\vth_0)\wedge (x-\vth_0))
  -F(X_t-\vth_0)F(x-\vth_0)}{f(X_t-\vth_0)}\d W_t\\
  &=&\frac{1}{\sqrt T}\left(\int_{0}^{X_T}H(z-\vth_0,x-\vth_0)\d
  z-\int_0^{X_0}H(z-\vth_0,x-\vth_0)\d z\right)\\
  &&-\frac{1}{\sqrt{T}}\int_{0}^{T}H(X_t-\vth_0,x-\vth_0)\d W_t.
\end{eqnarray*}
Using \eqref{extremeS} we have, for $x>A$,
\[1-F(x)=C\int_x^\infty{\rm exp}\left(2\int_0^yS_*(r)\d r\right)\d y
\leq C{\rm e}^{-2\gamma x},\] and
$$
  \frac{1-F(x)}{f(x)}
  \leq C\int_x^\infty{\rm e}^{-2\gamma(y-x)}\d y\leq C.
$$
For $x<-A$ we have $F(x)\leq C{\rm e}^{-2\gamma|x|}$ and we can write
\begin{eqnarray*}
  \frac{F(x)}{f(x)}=C\int_{-\infty}^x{\rm exp}(2\int_x^yS_*(r)\d r)\d
  y\leq C.
\end{eqnarray*}
These inequalities allow  us to deduce the following bounds
\begin{equation}
  \E_{\vth_0}H(\xi_{\vth_0}-\vth_0,x)^2=\E_0 H(\xi_0,x)^2<{\rm e}^{-\gamma |x|},
  \quad |x|>A.
\end{equation}
and
\begin{equation}
 \label{EDFextrem3}
  \E_{\vth_0}\left(\int_{0}^{\xi_{\vth_0}-\vth_0}H(z,x)\d z\right)^2
  =\E_0\left(\int_{0}^{\xi_0}H(z,x)\d z\right)^2\leq C{\rm e}^{-\gamma |x|},
  \quad |x|>A.
\end{equation}
Moreover
\begin{equation}
 \label{EDFextrem0}
  \int_{-\infty}^\infty\E_{0}\left(\int_{0}^{\xi_{0}}H(z,x)\d
  z\right)^2\d x\leq \infty.
\end{equation}
 Hence we get the asymptotic normality of $\eta_T^F(x)$:
\[
  \eta_T^F(x)\Longrightarrow \eta^F(x-\vth_0)\sim\mathcal {N}(0,4\E_0
  \left(H(\xi_0,x-\vth_0)\right)^2).
\]

As in Lemma \ref{LemConv} and Lemma \ref{LemIntConv}, if conditions $\mathcal {A}$ and $\mathcal {A}_0$ hold, we can show the
convergence of the vector $(\eta_T^F(x_1),...,\eta_T^F(x_k),\hat
u_T)$:
$$
  \mathcal{L}_{\vth_0}\left(\eta_T^F(x_1),...,\eta_T^F(x_k),\hat u_T\right)\Longrightarrow
  \mathcal{L}_{\vth_0}\left(\eta^F(x_1-\vth_0),...,\eta_T^F(x_k-\vth_0),\hat u\right)
$$
and the convergence of the integral:
\begin{eqnarray*}
  \mathcal{L}_{\vth_0}\big\{\int_{-\infty}^\infty
  \left(\eta_T^F(x)-\hat u_Tf(x-\vth_0)\right)^2\d
  x\big\}&\Longrightarrow&
  \mathcal{L}\left\{\int_{-\infty}^\infty
  \left(\eta^F(x)-\hat u f(x)\right)^2\d x\right\}.
\end{eqnarray*}
 We obtain
finally
\begin{eqnarray*}
  \Delta_T&=&T\int_{-\infty}^\infty(\hat F_T(x)-F(x,\hat\vth_T))^2\d x\\
  &=&\int_{-\infty}^\infty\big[\sqrt T(\hat F_T(x)-F(x-\vth_0))
  -\sqrt T(\hat \vth_T-\vth_0)F'(x-\tilde\vth_T)\big]^2\d
  x\\
  &=&\int_{-\infty}^\infty\big[\eta_T^F(x)-\hat u_Tf(x-\tilde\vth_T)\big]^2\d
  x\\
  &=&\int_{-\infty}^\infty\big[\eta_T^F(x)
  -\hat u_Tf(x-\vth_0)\big]^2\d x+o(1)\\
  &\Longrightarrow&\int_{-\infty}^\infty\big[\eta^F(x-\vth_0)
  -\hat u f(x-\vth_0)\big]^2\d x\\
  &=&\int_{-\infty}^\infty\left(\eta^F(y)-\hat u f(y)\right)^2\d
  y=\delta.
\end{eqnarray*}
So that the limit of the statistic $\Delta $ does not depend on
$\vth_0$, and the test $\Psi_T=\1_{\{\Delta_T\geq c_\vep\}}$ with
$c_\vep$ the solution of
$$
  \P\left(\Delta\geq c_\vep\right)=\vep
$$
belongs to $\mathcal {K}_\vep$.\\

{\bf Remark.}  It can be shown that in the case of Kolmogorov-Smirnov tests
$$
\varphi _T=\1_{\left\{\omega  _T>p_\varepsilon \right\}},\qquad
\quad \Phi _T=\1_{\left\{ \Omega _T>q_\varepsilon \right\}}
$$
where
$$
\omega  _T=\sup_x\left|\hat f_T\left(x\right)-f\left(x-\hat\vartheta\right)\right|\sqrt{T},\qquad \Omega  _T=\sup_x\left|\hat F_T\left(x\right)-F\left(x-\hat\vartheta\right)\right|\sqrt{T}
$$
the limit distributions of these statistics (under hypothesis) do not depend
on $\vartheta  $. The proofs can be done following the same lines as in
Kutoyants \cite{Kuto} and Negri \cite{Negri98} respectively.

\section{Consistency}

In this section we discuss the consistency of the proposed tests. We study the tests statistics under the alternative hypothesis that is defined as
$$
\mathcal {H}_1:\ S(\cdot)\not\in \overline{\mathcal {S}(\Theta)},
$$
where  $\overline{\mathcal {S}(\Theta)}=\left\{S\left(x-\vartheta
\right),\vartheta \in \left[\alpha ,\beta \right] \right\} $.

Under this hypothesis we have:
\begin{prop}
  Let all drift coefficients under alternative
  satisfy the conditions  $\mathcal {E}\mathcal {S},\ \mathcal {A}_0,$ and $\mathcal
  {A}$, then for any $S(\cdot)\not\in\overline{\mathcal {S}(\Theta)}$  we have
\end{prop}
  $$
   \P_S\left(\delta_T>d_\vep\right)\longrightarrow1,
  $$
  and
  $$
  \P_S\left(\Delta_T>c_\vep\right)\longrightarrow1.
  $$

\noindent\textbf{Proof.} Remember that under hypothesis $\mathcal
{H}_1$, the MLE $\hat\vth_T$ converges to the point which minimize
the distance
$$
D(\vth)=\E_S\left(S_*(\xi-\vth)-S(\xi)\right)^2,
$$
where $\xi$ is the random variable of invariant density $f_S(x)$
(See \cite{Kuto}, Proposition 2.36):
$$
\hat\vth_T \longrightarrow  \hat\vth_0={\rm
arg}\inf\limits_{\vth\in\Theta}D(\vth).
$$
In addition, denoted with $\|\cdot\|$ the norm in $L^2$, we have
\begin{eqnarray*}
&&\P_S\left(\delta_T>d_\vep\right)=\P_S\left
(\left\|\hat f_T(\cdot)-f(\cdot,\hat\vth_T)\right\|^2>d_\vep\right)\\
&&\quad\geq\P_S\left(\left\|f_S(x)-f(x-\hat\vth_T)\right\|^2-\left\|\hat
f_T(x)-f_S(x)\right\|^2>d_\vep\right).
\end{eqnarray*}
We can deduce
\begin{eqnarray*}
  &&\left\|f_S(x)-f(x-\hat\vth_T)\right\|^2=T\int_{-\infty}^\infty
  \left(f_S(x)-f(x-\hat\vth_T)\right)^2\d x\\
  &&\quad= T\int_{-\infty}^\infty
  \left(f_S(x)-f(x-\hat\vth_0)+o(1)\right)^2\d x\\
  &&\quad = (C+o(1))T\longrightarrow\infty,\quad {\rm as\ }T\longrightarrow\infty.
\end{eqnarray*}
Moreover
\begin{eqnarray*}
  &&\E_S\left(\left\|\hat
f_T(x)-f_S(x)\right\|^2\right)=\E_S\left(T\int_{-\infty}^\infty
  \left(\hat f_T(x)-f_S(x)\right)^2\d x\right)\\
  &&\quad\leq C\int_{-\infty}^\infty\E_S(\eta_T(x)^2)\d x
  \leq C\int_{-\infty}^\infty{\rm e}^{-2\gamma \left|x\right|}\d x<\infty.
\end{eqnarray*}
And finally we have the result for $\delta_T$:
  $$
   \P_S\left(\delta_T>d_\vep\right)\geq\P_S\left(\left\|f_S(x)-f(x-\hat\vth_T)\right\|^2-\left\|\hat
f_T(x)-f_S(x)\right\|^2>d_\vep\right)\longrightarrow1.
  $$
A similar result can be obtained for $\Delta_T$.

\section{Numerical Example}

 We
consider the Ornstein-Uhlenbeck process. Remind that the tests for O-U process
were studied in \cite{Kutoyants2012} as well. Suppose that the
observed process under the null hypothesis is

\[
\d X_t=-(X_t-\vth_0)\d t+\d W_t,\quad X_0,\ 0\leq t\leq T.
\]
The invariant density is $f(x-\vth_0)$, where
$f(x)= \pi^{-1/2}{\rm e}^{-x^2}$.

The log-likelihood ratio is

\[L(X^T,\vth)=-\int_0^T(X_t-\vth)\d X_t-\frac{1}{2}\int_0^T(X_t-\vth)^2\d t,\]
so that the MLE $\hat \vth_T$  can be calculated as

\[\hat \vth_T=\frac{1}{T}\int_0^TX_t\d t+\frac{X_T-X_0}{T}.\]
The Fisher information in this case equals to 1, and the LTE is

\[
\hat f_T(x)=\frac{1}{T}(|X_T-x|-|X_0-x|)-\frac{1}{T}\int_0^T{\rm
sgn}(X_t-x)\d X_t.
\]

The conditions $\mathcal {A}_0$ and $\mathcal {A}$ are
fulfilled, then the statistic is convergent:
\[
\delta_T=\int_{-\infty}^\infty\left(\hat
f_T(x)-f(x-\hat\vth_T)\right)^2\d x\Longrightarrow
\delta=\int_{-\infty}^\infty\zeta_1(x)^2\d x,
\]
where the limit process $\zeta_1(x)=\eta(x)-\hat u f'(x)$ can be
written as
$$
\zeta_1(x)= \int_{-\infty}^\infty
\left(2f(x)\frac{\1_{\{y>x\}}-F(y)}{\sqrt{f(y)}}+f'(x)\sqrt{f(y)}\right)\d
W(y).\\
$$

We have a similar result for the test based on the EDF:
$$
\Delta_T=\int_{-\infty}^\infty\left(\hat
F_T(x)-F(x-\hat\vth_T)\right)^2\d x \Longrightarrow
\Delta=\int_{-\infty}^\infty\left(\zeta_2(x)\right)^2\d x,\\
$$
where the limit process can be written as
$$
\zeta_2(x)=\int_{-\infty}^\infty \left(2\frac{F(y\wedge
x)-F(y)F(x)}{\sqrt{f(y)}}+f(x)\sqrt{f(y)}\right)\d W(y).
$$

\bigskip

\begin{figure}
\centering
\includegraphics[width=5.77in,height=3.00in]{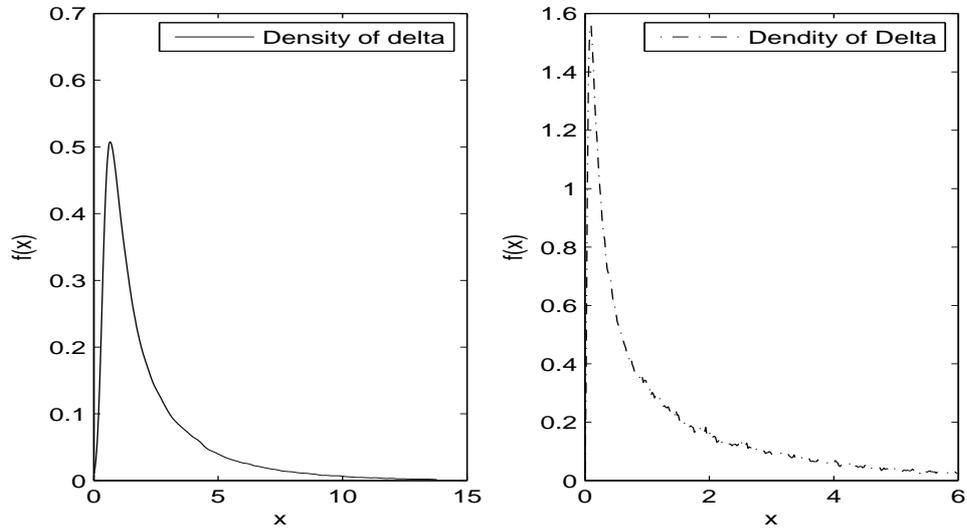}
\caption{Density of the statistics. On the left the density of $\delta$, on the right the density of $\Delta$} \label{density}
\end{figure}

We simulate $10^5$ trajectories of $\delta$ (resp. $\Delta$) and
calculate the empirical $1-\vep$ quantiles of $\delta$ (resp.
$\Delta$).
We obtain the simulated density for $\delta$ and $\Delta$ that are showed in  Graphic
\ref{density}. The values of the  thresholds $d_\vep $ for different $\vep$ are showed in
Graphic \ref{d_vep}.

\begin{figure}
\centering
\includegraphics[width=5.77in,height=3.00in]{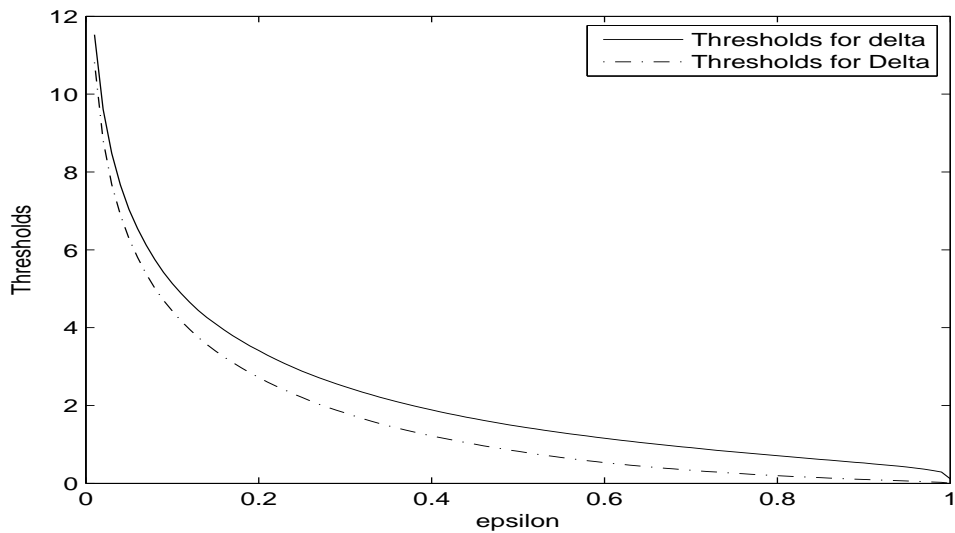}
\caption{Threshold for different $\vep$. The solid line represents  the values
for $\delta$, the dotted line represents the values for $\Delta$} \label{d_vep}
\end{figure}

\subsection*{Acknowledgments}
The authors are grateful to Y. Kutoyants for his suggestions and interesting discussions.

\end{document}